\documentclass{amsart}
\usepackage{amsmath,amssymb}
\usepackage{yhmath}

\newcommand{\bdis}{\begin{displaymath}}
\newcommand{\edis}{\end{displaymath}}
\newcommand{\be}{\begin{equation}}
\newcommand{\ee}{\end{equation}}
\newcommand{\mbb}{\mathbb}
\newcommand{\mcal}{\mathcal}

\newcommand{\vp}{\varphi} 
\newcommand{\vth}{\vartheta}

\newcommand{\zf}{\zeta\left(\frac{1}{2}+it\right)}

\newcommand{\zfvpr}{\zeta\left(\frac{1}{2}+i\vp_1^r(t)\right)}

\newcommand{\FR}{\frac{x^n+y^n}{z^n}}

\DeclareMathOperator{\im}{Im}

\theoremstyle{definition}

\theoremstyle{remark}
\newtheorem{remark}[]{Remark}

\newtheorem*{mydef11}{{\bf Theorem 1}}

\newtheorem*{mydef12}{{\bf Theorem 2}}

\newtheorem*{mydef13}{{\bf Theorem 3}}

\newtheorem*{mydef14}{{\bf Theorem 4}}

\newtheorem*{mydef15}{{\bf Theorem 5}}

\newtheorem*{mydef41}{{\bf Corollary 1}}

\newtheorem*{mydef42}{{\bf Corollary 2}}

\newtheorem*{mydef43}{{\bf Corollary 3}}

\newtheorem*{mydef51}{{\bf Lemma 1}}

\newtheorem*{mydef81}{{\bf Property 1}}

\newtheorem*{mydef82}{{\bf Property 2}}

\numberwithin{equation}{section}

\begin{document}

\title[Jacob's ladders, our asymptotic formulae (1981)  \dots]{Jacob's ladders, our asymptotic formulae (1981) connected with the equation $Z'(t)=0$ as a base for new $\zeta$-equivalents of the Fermat-Wiles theorem and some other results}

\author{Jan Moser}

\address{Department of Mathematical Analysis and Numerical Mathematics, Comenius University, Mlynska Dolina M105, 842 48 Bratislava, SLOVAKIA}

\email{jan.mozer@fmph.uniba.sk}


\begin{abstract}
In this paper we obtain new functionals and corresponding $\zeta$-equivalents of the Fermat-Wiles theorem. These are generated by our asymptotic formulae (1981) for the function $Z'(t)$ which we are now using to prove essential influence of the Lindel\" of hypothesis on the distribution of the roots of odd order of the equation $Z'(t)=0$. 
\end{abstract}
\maketitle

\section{Introduction} 

\subsection{} 

Let us remind that Riemann has defined also the following real-valued function 
\be \label{1.1} 
Z(t)=e^{i\vth(t)}\zf, 
\ee  
where 
\be \label{1.2} 
\begin{split}
& \vth(t)=-\frac t2\ln\pi+\im\Gamma\left(\frac 14+i\frac t2\right)= \\ 
& \frac t2\ln\frac{t}{2\pi}-\frac t2-\frac{\pi}{8}+\mcal{O}\left(\frac 1t\right) 
\end{split}
\ee 
(see \cite{11}, (35), (44), (62), comp. \cite{12}, p. 98). 

Next, let 
\be \label{1.3} 
S(a,b)=\sum_{0<a\leq n<b\leq 2a}n^{it},\ b\leq \sqrt{\frac{t}{2\pi}}
\ee 
be the elementary trigonometric sum. In our paper \cite{4}, (2), (16) we have obtained a completely new asymptotic formulae: if 
\be \label{1.4} 
|S(a,b)|<A(\Delta)\sqrt{a}t^\Delta,\ 0<\Delta\leq \frac 16, 
\ee  
then 
\be \label{1.5} 
\begin{split} 
& \sum_{T\leq\bar{t}_{2\nu}\leq T+H}Z'(\bar{t}_{2\nu})=-\frac{1}{4\pi}H\ln^2\frac{T}{2\pi}+\mcal{O}(T^{\Delta}\ln^2T), \\ 
& \sum_{T\leq\bar{t}_{2\nu+1}\leq T+H}Z'(\bar{t}_{2\nu+1})=-\frac{1}{4\pi}H\ln^2\frac{T}{2\pi}+\mcal{O}(T^{\Delta}\ln^2T), 
\end{split} 
\ee 
where the sequence 
\be \label{1.6} 
\{\bar{t}_\nu\}_{\nu=1}^\infty 
\ee 
is defined by the condition 
\be \label{1.7} 
\vth(\bar{t}_\nu)=\pi\nu+\frac{\pi}{2}, 
\ee   
and further we require also that 
\be \label{1.8} 
H=T^\Delta\psi(T), 
\ee  
where $\psi(T)$ is a function that is increasing with arbitrary slow growth to $+\infty$. 

\begin{remark}
In this paper, it is sufficient to put $\psi(T)=\ln T$, i.e. we have 
\be \label{1.9} 
H=T^\Delta\ln T. 
\ee 
\end{remark} 

The proof of our formulae (\ref{1.5}) is based on the following Formula 1 (see \cite{4}, (10), and the proof: (17) -- (39)) 
\be \label{1.10} 
Z'(t)=-2\sum_{n\leq\sqrt{\frac{t}{2\pi}}}\frac{1}{\sqrt{n}}(\vth'-\ln n)\sin(\vth-t\ln n)+\mcal{O}(t^{-1/4}\ln t), 
\ee  
that we have transformed into the Formula 2 (see \cite{4}, (11)) 
\be \label{1.11} 
\begin{split} 
& Z'(t)=-2\sum_{n<P_0}\frac{1}{\sqrt{n}}\ln\frac{P_0}{n}\sin(\vth-t\ln n)+\mcal{O}(T^{-1/4}\ln T), \\ 
& t\in [T,T+H],\ H\in (0,\sqrt[4]{T}],\ P_0=\sqrt{\frac{T}{2\pi}}. 
\end{split} 
\ee 

\begin{remark}
The factor in the amplitude of Riemann's oscillators  
\be \label{1.12} 
\ln\frac{P_0}{n},\ n<P_0 
\ee 
is the matter of prime importance in our formula (\ref{1.11}) since it enabled us to remove small denominators, comp. \cite{4}, (57) -- (64), in the course of estimating the corresponding sums. 
\end{remark} 

\subsection{} 

In this paper we prove that our asymptotic formulae (\ref{1.5}), connected with estimates of elementary trigonometric sums\footnote{Comp. (\ref{1.3}) and (\ref{1.4}).}, also generates, after 14 years, new functionals together with corresponding new $\zeta$-equivalents of the Fermat-Wiles theorem. 

For example, we have obtained the following functional 
\be \label{1.13} 
\begin{split}
& \lim_{\tau\to\infty}\frac{1}{\tau} 
\left\{
\sum_{\bar{t}_{2\nu+1}\geq (x\tau)^{1/\Delta}}^{\bar{t}_{2\nu+1}\leq (x\tau)^{1/\Delta}+H((x\tau)^{1/\Delta})}Z'(\bar{\tau}_{2\nu+1})
\right\}\times \\ 
& \left\{
\int_{[(x\tau)^{1/\Delta}]^3}^{[(x\tau)^{1/\Delta}+2\bar{l}]^3}\prod_{r=0}^2Z^2(\vp_1^r(t)){\rm d}t
\right\}^{-1}=x
\end{split}
\ee 
for every fixed $x>0$ and every admissible $\Delta\in (0,\frac 16]$, where 
\be \label{1.14} 
H(T)=T^\Delta\ln T,\ 2\bar{l}=\frac{1}{4\pi},\ [G]^3=\vp_1^{-3}(G). 
\ee 

\begin{remark}
Admissible $\Delta\in (0,1/6]$ denotes those $\Delta$ for which the estimate (\ref{1.4}) holds true. 
\end{remark} 

In the special case 
\be \label{1.15} 
x\to\FR,\ x,y,z,n\in\mbb{N},\ n\geq 3 
\ee  
of the Fermat's rationals we have the following result: the $\zeta$-condition 
\be \label{1.16} 
\begin{split}
& \lim_{\tau\to\infty}\frac{1}{\tau} 
\left\{
\sum_{\bar{t}_{2\nu+1}\geq (\FR\tau)^{1/\Delta}}^{\bar{t}_{2\nu+1}\leq (\FR\tau)^{1/\Delta}+H((\FR\tau)^{1/\Delta})}Z'(\bar{\tau}_{2\nu+1})
\right\}\times \\ 
& \left\{
\int_{[(\FR\tau)^{1/\Delta}]^3}^{[(\FR\tau)^{1/\Delta}+2\bar{l}]^3}\prod_{r=0}^2Z^2(\vp_1^r(t)){\rm d}t
\right\}^{-1}\not=1
\end{split}
\ee 
on the set of all Fermat's rationals and for every fixed and admissible $\Delta\in(0,1/6]$ expresses new $\zeta$-equivalent of the Fermat-Wiles theorem. 

Next, we have obtained the asymptotic formula of such a kind: 
\be \label{1.17} 
\begin{split}
& \sum_{\bar{t}_{2\nu+1}\geq (\tau)^{1/\Delta}}^{\bar{t}_{2\nu+1}\leq (\tau)^{1/\Delta}+H((\tau)^{1/\Delta})}Z'(\bar{\tau}_{2\nu+1})\sim \\ 
& \int_{\tau}^{\overset{1}{\tau}}Z^2(t){\rm d}d\times \int_{[(\tau)^{1/\Delta}]^3}^{[(\tau)^{1/\Delta}+2\tilde{l}]^2}\prod_{r=0}^2Z^2(\vp_1^r(t)){\rm d}t, \\ 
& 2\tilde{l}=\frac{1}{4\pi(1-c)}
\end{split}
\ee 
as $\tau\to\infty$ and for every fixed admissible $\Delta\in (0,1/6]$. 

\begin{remark}
\emph{An asymptotic equilibrium} of the three mathematical objects 
\be \label{1.18} 
\begin{split} 
& \left\{
\sum_{\bar{t}_{2\nu+1}\geq (\tau)^{1/\Delta}}^{\bar{t}_{2\nu+1}\leq (\tau)^{1/\Delta}+H((\tau)^{1/\Delta})}Z'(\bar{\tau}_{2\nu+1}),\ \int_{\tau}^{\overset{1}{\tau}}Z^2(t){\rm d}t,\right. \\  
& \left. \int_{[(\tau)^{1/\Delta}]^3}^{[(\tau)^{1/\Delta}+2\tilde{l}]^2}\prod_{r=0}^2Z^2(\vp_1^r(t)){\rm d}t 
\right\} 
\end{split} 
\ee 
is expressed by the formula (\ref{1.17}). 
\end{remark} 

There is also the following meaning of the formula (\ref{1.17}). Let us put 
\be \label{1.19} 
C_1=\int_{\tau}^{\overset{1}{\tau}}Z^2(t){\rm d}t,\ C_2=\int_{[(\tau)^{1/\Delta}]^3}^{[(\tau)^{1/\Delta}+2\tilde{l}]^2}\prod_{r=0}^2Z^2(\vp_1^r(t)){\rm d}t, 
\ee  
where $C_1+C_2=C$ is the hypotenuse of the corresponding right-angled triangle. Then, by Euclid's theorem 
\bdis 
V^2=C_1C_2
\edis  
for the altitude of this triangle. 

\begin{remark}
Consequently, we have the following asymptotic formula by (\ref{1.17}). 
\be \label{1.20} 
V=\left\{
\sum_{\bar{t}_{2\nu+1}\geq (\tau)^{1/\Delta}}^{\bar{t}_{2\nu+1}\leq (\tau)^{1/\Delta}+H((\tau)^{1/\Delta})}Z'(\bar{\tau}_{2\nu+1})
\right\}^{1/2},\ \tau\to\infty
\ee 
for this Euclidean altitude. Of course, the formula (\ref{1.20}) represents also the asymptotic value of the geometric mean of the two integrals in (\ref{1.19}). 
\end{remark} 

\section{On I. M. Vinogradov' scepticism on possibilities of the method of trigonometric sums and his analogue for elementary trigonometric sums} 

\subsection{} 

I. M. Vinogradov has analysed in the Introduction to his monograph \cite{13} the possibilities of the method of trigonometric sums (H. Weyl's sums) in the problem of estimating the reminder term $R(N)$ in the asymptotic formula\footnote{See \cite{13}, p. 13.} 
\be \label{2.1} 
\pi(N)-\int_2^{N}\frac{{\rm d}x}{\ln x}=R(N), 
\ee  
where we have: 
\begin{itemize}
	\item[(a)] The original Vall\` e-Poussin's result 
	\be \label{2.2} 
	R(N)=\mcal{O}(Ne^{-c\sqrt{\ln N}}), 
	\ee  
	\item[(b)] H. Weyl's result 
	\be \label{2.3} 
	R(N)=\mcal{O}(Ne^{-c_1\sqrt{\ln N\ln\ln N}}), 
	\ee  
	\item[(c)] 
	Vinogradov' result 
	\be \label{2.4} 
	\begin{split} 
	& R(N)=\mcal{O}(Ne^{-c_2\lambda(N)\ln^{0.6}N}), \\ 
	& \lambda(N)=(\ln\ln N)^{-0.2}. 
	\end{split} 
	\ee 
\end{itemize} 

From the Riemann hypothesis, however, the next estimate follows 
\be \label{2.5} 
R(N)=\mcal{O}(\sqrt{N}\ln N), 
\ee  
and this one is extremely distanced from the estimate (\ref{2.4}). In connection with this fact, I. M. Vinogradov made the remark: \emph{Obviously, it is very hard to make an essential progress in solution of the problem to find the order of the $N$-th term willing to find} 
\be \label{2.6} 
R(N)=\mcal{O}(N^{1-\gamma}),\ \gamma=0.000001
\ee 
\emph{by making use of only some improvements of the H. Weyl's estimates and without making use of further important progresses in the theory of zeta-function}. 

\subsection{} 

The main Theorem of our paper \cite{4} follows from the formulae (\ref{1.5}): If 
\be \label{2.7} 
|S(a,b)|< A(\Delta)\sqrt{a}t^\Delta,\ 0<\Delta\leq \frac 16, 
\ee  
then there is a root of odd order of the equation 
\be \label{2.8} 
Z'(t)=0
\ee 
in the interval\footnote{Comp. (\ref{1.8}) and (\ref{1.9}).} 
\be \label{2.9} 
(T,T+T^\Delta\ln T). 
\ee  
For example, in the case of Kolesnik's exponent, see \cite{3}, 
\be \label{2.10} 
\Delta=\frac{35}{216}+\epsilon, 
\ee  
we have the following Corollary 1: The interval 
\be \label{2.11} 
(T,T+T^{\frac{35}{216}+\epsilon})
\ee  
contains the root of odd order of the equation (\ref{2.8}). Here $0<\epsilon$ be arbitrarily small. 

\begin{remark}
The Kolesnik's value $\frac{35}{216}$ is only $2.8\%$ improvement of the exponent $\frac{1}{6}$. 
\end{remark}  

On the Lindel\" of hypothesis we have\footnote{Comp. \cite{2}, p. 89.} that 
\be \label{2.12} 
|S(a,b)|<A(\epsilon)\sqrt{a}t^\epsilon, 
\ee  
for every sufficiently small and fixed $\epsilon>0$, i.e. we have the following Corollary 2: On the Lindel\" of hypothesis the interval 
\be \label{2.13} 
(T,T+T^\epsilon\ln T)
\ee 
contains a root of odd order of the equation (\ref{2.8}). 

\begin{remark}
Thus, the Lindel\" of hypothesis gives an $100\%$ improvement of the Kolesnik's exponent in (\ref{2.11}). Simultaneously, we have an $100\%$ improvement of all possible next improvements of the Kolesnik's exponent for the elementary trigonometric sums. 
\end{remark} 

\subsection{} 

In connection with our above-mentioned Theorem of the paper \cite{4} we can give the following analogue of the I. M. Vinogradov' scepticism. 

\begin{remark}
It is very hard to find, for example, the following estimate 
\be \label{2.14} 
|S(a,b)|<A_1\sqrt{a}t^{0.000001},\ t\to\infty
\ee 
by making use of only some improvements of the Weyl's estimates. 
\end{remark}

\section{Jacob's ladders: notions and basic geometrical properties}  

\subsection{}

In this paper we use the following notions of our works \cite{5} -- \cite{9}: 
\begin{itemize}
\item[{\tt (a)}] Jacob's ladder $\vp_1(T)$, 
\item[{\tt (b)}] direct iterations of Jacob's ladders 
\bdis 
\begin{split}
	& \vp_1^0(t)=t,\ \vp_1^1(t)=\vp_1(t),\ \vp_1^2(t)=\vp_1(\vp_1(t)),\dots , \\ 
	& \vp_1^k(t)=\vp_1(\vp_1^{k-1}(t))
\end{split}
\edis 
for every fixed natural number $k$, 
\item[{\tt (c)}] reverse iterations of Jacob's ladders 
\be \label{3.1}  
\begin{split}
	& \vp_1^{-1}(T)=\overset{1}{T},\ \vp_1^{-2}(T)=\vp_1^{-1}(\overset{1}{T})=\overset{2}{T},\dots, \\ 
	& \vp_1^{-r}(T)=\vp_1^{-1}(\overset{r-1}{T})=\overset{r}{T},\ r=1,\dots,k, 
\end{split} 
\ee   
where, for example, 
\be \label{3.2} 
\vp_1(\overset{r}{T})=\overset{r-1}{T}
\ee  
for every fixed $k\in\mbb{N}$ and every sufficiently big $T>0$. We also use the properties of the reverse iterations listed below.  
\be \label{3.3}
\overset{r}{T}-\overset{r-1}{T}\sim(1-c)\pi(\overset{r}{T});\ \pi(\overset{r}{T})\sim\frac{\overset{r}{T}}{\ln \overset{r}{T}},\ r=1,\dots,k,\ T\to\infty,  
\ee 
\be \label{3.4} 
\overset{0}{T}=T<\overset{1}{T}(T)<\overset{2}{T}(T)<\dots<\overset{k}{T}(T), 
\ee 
and 
\be \label{3.5} 
T\sim \overset{1}{T}\sim \overset{2}{T}\sim \dots\sim \overset{k}{T},\ T\to\infty.   
\ee  
\end{itemize} 

\begin{remark}
	The asymptotic behaviour of the points 
	\bdis 
	\{T,\overset{1}{T},\dots,\overset{k}{T}\}
	\edis  
	is as follows: at $T\to\infty$ these points recede unboundedly each from other and all together are receding to infinity. Hence, the set of these points behaves at $T\to\infty$ as one-dimensional Friedmann-Hubble expanding Universe. 
\end{remark}  

\subsection{} 

Let us remind that we have proved\footnote{See \cite{9}, (3.4).} the existence of almost linear increments 
\be \label{3.6} 
\begin{split}
& \int_{\overset{r-1}{T}}^{\overset{r}{T}}\left|\zf\right|^2{\rm d}t\sim (1-c)\overset{r-1}{T}, \\ 
& r=1,\dots,k,\ T\to\infty,\ \overset{r}{T}=\overset{r}{T}(T)=\vp_1^{-r}(T)
\end{split} 
\ee 
for the Hardy-Littlewood integral (1918) 
\be \label{3.7} 
J(T)=\int_0^T\left|\zf\right|^2{\rm d}t. 
\ee  

For completeness, we give here some basic geometrical properties related to Jacob's ladders. These are generated by the sequence 
\be \label{3.8} 
T\to \left\{\overset{r}{T}(T)\right\}_{r=1}^k
\ee 
of reverse iterations of the Jacob's ladders for every sufficiently big $T>0$ and every fixed $k\in\mbb{N}$. 

\begin{mydef81}
The sequence (\ref{3.8}) defines a partition of the segment $[T,\overset{k}{T}]$ as follows 
\be \label{3.9} 
|[T,\overset{k}{T}]|=\sum_{r=1}^k|[\overset{r-1}{T},\overset{r}{T}]|
\ee 
on the asymptotically equidistant parts 
\be \label{3.10} 
\begin{split}
& \overset{r}{T}-\overset{r-1}{T}\sim \overset{r+1}{T}-\overset{r}{T}, \\ 
& r=1,\dots,k-1,\ T\to\infty. 
\end{split}
\ee 
\end{mydef81} 

\begin{mydef82}
Simultaneously with the Property 1, the sequence (\ref{3.8}) defines the partition of the integral 
\be \label{3.11} 
\int_T^{\overset{k}{T}}\left|\zf\right|^2{\rm d}t
\ee 
into the parts 
\be \label{3.12} 
\int_T^{\overset{k}{T}}\left|\zf\right|^2{\rm d}t=\sum_{r=1}^k\int_{\overset{r-1}{T}}^{\overset{r}{T}}\left|\zf\right|^2{\rm d}t, 
\ee 
that are asymptotically equal 
\be \label{3.13} 
\int_{\overset{r-1}{T}}^{\overset{r}{T}}\left|\zf\right|^2{\rm d}t\sim \int_{\overset{r}{T}}^{\overset{r+1}{T}}\left|\zf\right|^2{\rm d}t,\ T\to\infty. 
\ee 
\end{mydef82} 

It is clear, that (\ref{3.10}) follows from (\ref{3.3}) and (\ref{3.5}) since 
\be \label{3.14} 
\overset{r}{T}-\overset{r-1}{T}\sim (1-c)\frac{\overset{r}{T}}{\ln \overset{r}{T}}\sim (1-c)\frac{T}{\ln T},\ r=1,\dots,k, 
\ee  
while our eq. (\ref{3.13}) follows from (\ref{3.6}) and (\ref{3.5}).  

\section{The first new $\zeta$-functional and corresponding $\zeta$-equivalent of the Fermat-Wiles theorem generated by our formulae (\ref{1.5})} 

\subsection{} 

We choose the second formula in (\ref{1.5}) as the basic formula\footnote{See also (\ref{1.9}).}, and namely in the form 
\be \label{4.1} 
\begin{split}
& \sum_{\bar{t}_{2\nu+1}\geq T}^{\leq T+H(T)}Z'(\bar{t}_{2\nu+1})=\frac{1}{4\pi}T^\Delta\ln^3T+\mcal{O}(T^\Delta\ln^2T)= \\ 
& \left\{1+\mcal{O}\left(\frac{1}{\ln T}\right)\right\}\frac{1}{4\pi}T^\Delta\ln^3T. 
\end{split}
\ee 
Next, we use our formula\footnote{See \cite{10}, (3.18); $f_m\equiv 1$.}  
\be \label{4.2} 
\begin{split} 
& \int_{\overset{k}{T}}^{\overset{k}{\wideparen{T+2l}}}\prod_{r=0}^{k-1}\left|\zfvpr\right|^2{\rm d}r= \\ 
& \left\{1+\mcal{O}\left(\frac{\ln\ln T}{\ln T}\right)\right\}2l\ln^kT
\end{split}
\ee 
for every fixed $l>0$, $k\in\mbb{N}$ in the case $k=3$, and in the following form\footnote{Comp. (\ref{1.1}).} 
\be \label{4.3} 
2l\ln^3T=\left\{1+\mcal{O}\left(\frac{\ln\ln T}{\ln T}\right)\right\}\int_{\overset{3}{T}}^{\overset{3}{\wideparen{T+2l}}}\prod_{r=0}^{2}Z^2(\vp_1^r(t)){\rm d}t. 
\ee 
Now, it follows from (\ref{4.1}) by (\ref{4.3}) 
\begin{mydef51}
\be \label{4.4} 
\begin{split}
& \sum_{\bar{t}_{2\nu+1}\geq T}^{\leq T+H(T)}Z'(\bar{t}_{2\nu+1})=\\ 
& \left\{1+\mcal{O}\left(\frac{\ln\ln T}{\ln T}\right)\right\}\frac{1}{8\pi l}T^\Delta\int_{\overset{3}{T}}^{\overset{3}{\wideparen{T+2l}}}\prod_{r=0}^{2}Z^2(\vp_1^r(t)){\rm d}t
\end{split}
\ee 
for every fixed $l>0$, and every fixed admissible $\Delta\in (0,1/6]$, comp. Remark 3. 
\end{mydef51} 

Next, we choose the value $l$ as 
\be \label{4.5} 
\bar{l}=\frac{1}{8\pi}, 
\ee  
and after this we use the substitution 
\be \label{4.6} 
T=(x\tau)^{\frac{1}{\Delta}},\ x>0 
\ee 
in the formula (\ref{4.4}), $l=\bar{l}$, that gives us the following functional. 

\begin{mydef11}
\be \label{4.7} 
\begin{split}
& \lim_{\tau\to\infty}\frac{1}{\tau}\left\{
\sum_{\bar{t}_{2\nu+1}\geq (x\tau)^{1/\Delta}}^{\leq (x\tau)^{1/\Delta}+H((x\tau)^{1/\Delta})}Z'(\bar{t}_{2\nu+1})
\right\}\times \\ 
& \left\{
\int_{[(x\tau)^{1/\Delta}]^3}^{[(x\tau)^{1/\Delta}+2\bar{l}]^3}\prod_{r=0}^{2}Z^2(\vp_1^r(t)){\rm d}t
\right\}^{-1}=x
\end{split}
\ee 
where 
\be \label{4.8} 
H(T)=T^\Delta\ln T,\ 2\bar{l}=\frac{1}{4\pi}
\ee 
for every fixed $x>0$ and admissible $\Delta\in (0,1/6]$.  
\end{mydef11}

In the special case of Fermat's rationals 
\be \label{4.9} 
x\to\FR,\ x,y,z,n\in\mbb{N},\ n\geq 3
\ee 
we obtain the following corollary. 

\begin{mydef41}
\be \label{4.10} 
\begin{split}
& \lim_{\tau\to\infty}\frac{1}{\tau}\left\{
\sum_{\bar{t}_{2\nu+1}\geq (\FR\tau)^{1/\Delta}}^{\leq (\FR\tau)^{1/\Delta}+H((\FR\tau)^{1/\Delta})}Z'(\bar{t}_{2\nu+1})
\right\}\times \\ 
& \left\{
\int_{[(\FR\tau)^{1/\Delta}]^3}^{(\FR\tau)^{1/\Delta}+2\bar{l}]^3}\prod_{r=0}^{2}Z^2(\vp_1^r(t)){\rm d}t
\right\}^{-1} = \FR
\end{split}
\ee 
for every fixed Fermat's rational and every fixed admissible $\Delta\in (0,1/6]$. 
\end{mydef41} 

Consequently, the next theorem follows from (\ref{4.10}). 

\begin{mydef12}
The $\zeta$-condition 
\be \label{4.11} 
\begin{split}
& \lim_{\tau\to\infty}\frac{1}{\tau}\left\{
\sum_{\bar{t}_{2\nu+1}\geq (\FR\tau)^{1/\Delta}}^{\leq (\FR\tau)^{1/\Delta}+H((\FR\tau)^{1/\Delta})}Z'(\bar{t}_{2\nu+1})
\right\}\times \\ 
& \left\{
\int_{[(\FR\tau)^{1/\Delta}]^3}^{(\FR\tau)^{1/\Delta}+2\bar{l}]^3}\prod_{r=0}^{2}Z^2(\vp_1^r(t)){\rm d}t
\right\}^{-1} \not=1 
\end{split}
\ee 
on the set of all Fermat's rationals and every fixed  admissible $\Delta\in (0,1/6]$ expresses the new $\zeta$-equivalent of the Fermat-Wiles theorem. 
\end{mydef12} 

\subsection{} 

Here we geve, for example, some concrete functionals contained in the formula (\ref{4.7}). 

\begin{itemize}
	\item[(A)] In the case $\Delta=\frac 16$ we obtain the next result. 
	\begin{mydef42}  
	\be \label{4.12} 
	\begin{split}
	& \lim_{\tau\to\infty}\frac{1}{\tau}\left\{
	\sum_{\bar{t}_{2\nu+1}\geq (x\tau)^{6}}^{\leq (x\tau)^{6}+H((x\tau)^{6})}Z'(\bar{t}_{2\nu+1})
	\right\}\times \\ 
	& \left\{
	\int_{[(x\tau)^{6}]^3}^{(x\tau)^{6}+2\bar{l}]^3}\prod_{r=0}^{2}Z^2(\vp_1^r(t)){\rm d}t
	\right\}^{-1}=x.
	\end{split}
	\ee 
	\end{mydef42} 
	\item[(B)] On the Lindel\" of hypothesis we can put, for example, $\epsilon=10^{-6}$ and we obtain the next result. 
	\begin{mydef43}
	\be\label{4.13} 
	\begin{split} 
	& \lim_{\tau\to\infty}\frac{1}{\tau}\left\{
	\sum_{\bar{t}_{2\nu+1}\geq (x\tau)^{10^6}}^{\leq (x\tau)^{10^6}+H((x\tau)^{10^6})}Z'(\bar{t}_{2\nu+1})
	\right\}\times \\ 
	& \left\{
	\int_{[(x\tau)^{10^6}]^3}^{(x\tau)^{10^6}+2\bar{l}]^3}\prod_{r=0}^{2}Z^2(\vp_1^r(t)){\rm d}t
	\right\}^{-1}=x.
	\end{split} 
	\ee 
	\end{mydef43}
\end{itemize} 

\section{On \emph{EQUILIBRIUM STATE} generated by the formula (\ref{4.4})} 

\subsection{} 

We use the following two operations on the formula (\ref{4.4}): First, we use the substitution 
\be \label{5.1} 
T=\tau^{1/\Delta}, 
\ee  
that gives the result 
\be \label{5.2} 
\begin{split}
& \sum_{\bar{t}_{2\nu+1}\geq \tau^{1/\Delta}}^{\leq \tau^{1/\Delta}+H(\tau^{1/\Delta})}Z'(\bar{t}_{2\nu+1})=\\ 
& \left\{1+\mcal{O}\left(\frac{\ln\ln \tau}{\ln \tau}\right)\right\}\frac{1}{8\pi l}\tau\int_{[\tau^{1/\Delta}]^3}^{[\tau^{1/\Delta}+2l]^3}\prod_{r=0}^{2}Z^2(\vp_1^r(t)){\rm d}t, 
\end{split}
\ee  
and next, we use our almost linear formula\footnote{See \cite{9}, (3.4), (3.6) and $r=1$.} in the form 
\be \label{5.3} 
(1-c)\tau=\{1+\mcal{O}(\tau^{-2/3+\delta})\}\int_\tau^{\overset{1}{\tau}}Z^2(t){\rm d}t 
\ee  
($0<\delta$ being sufficiently small) that gives the result 
\be \label{5.4} 
\begin{split}
& \sum_{\bar{t}_{2\nu+1}\geq \tau^{1/\Delta}}^{\leq \tau^{1/\Delta}+H(\tau^{1/\Delta})}Z'(\bar{t}_{2\nu+1})=\\ 
& \left\{1+\mcal{O}\left(\frac{\ln\ln \tau}{\ln \tau}\right)\right\}\frac{1}{8\pi (1-c)l}\int_\tau^{\overset{1}{\tau}}Z^2(t){\rm d}t\times \\ 
& \int_{[\tau^{1/\Delta}]^3}^{[\tau^{1/\Delta}+2l]^3}\prod_{r=0}^{2}Z^2(\vp_1^r(t)){\rm d}t,\ \tau\to\infty. 
\end{split}
\ee 
Now, in the case 
\be \label{5.5} 
2\tilde{l}=\frac{1}{4\pi(1-c)}, 
\ee  
the next theorem follows from (\ref{5.4}). 

\begin{mydef13}
\be \label{5.6} 
\begin{split}
& \sum_{\bar{t}_{2\nu+1}\geq \tau^{1/\Delta}}^{\leq \tau^{1/\Delta}+H(\tau^{1/\Delta})}Z'(\bar{t}_{2\nu+1})\sim \\ 
& \int_\tau^{\overset{1}{\tau}}Z^2(t){\rm d}t\times \int_{[\tau^{1/\Delta}]^3}^{[\tau^{1/\Delta}+2l]^3}\prod_{r=0}^{2}Z^2(\vp_1^r(t)){\rm d}t, 
\end{split}
\ee 
as $\tau\to\infty$ for every fixed and admissible $\Delta\in(0,1/6]$. 
\end{mydef13} 

\begin{remark}
Some interpretations of the formula (\ref{5.6}) are contained in the Remarks 4 and 5, comp. (\ref{1.17}). 
\end{remark} 

\section{Results that are corresponding with the partition $3=2+1$ in the formula (\ref{4.1})} 

\subsection{} 

Now, we write the formula (\ref{4.1}) in the form 
\be \label{6.1} 
\begin{split}
& \sum_{\bar{t}_{2\nu+1}\geq T}^{\leq T+H(T)}Z'(\bar{t}_{2\nu+1})= \\ 
& \left\{1+\mcal{O}\left(\frac{1}{\ln T}\right)\right\}\frac{1}{4\pi}T^\Delta\ln^2T\ln T, 
\end{split}
\ee  
and next we use the following formulae\footnote{See (\ref{4.2}).} 
\be \label{6.2} 
\begin{split}
& 2l_1\ln T=\left\{1+\mcal{O}\left(\frac{\ln\ln T}{\ln T}\right)\right\}\int_{\overset{1}{T}}^{\overset{1}{\wideparen{T+2l_1}}}Z^2(t){\rm d}t, \\ 
& 2l_2\ln^2T=\left\{1+\mcal{O}\left(\frac{\ln\ln T}{\ln T}\right)\right\}\int_{\overset{2}{T}}^{\overset{2}{\wideparen{T+2l_2}}}\prod_{r=0}^1Z^2(\vp_1^r(t)){\rm d}t, 
\end{split}
\ee 
that gives the following result 
\be \label{6.3} 
\begin{split}
& \sum_{\bar{t}_{2\nu+1}\geq T}^{\leq T+H(T)}Z'(\bar{t}_{2\nu+1})= \\ 
& \left\{1+\mcal{O}\left(\frac{\ln\ln T}{\ln T}\right)\right\}\frac{T^\Delta}{16\pi l_1l_2}\times \int_{\overset{1}{T}}^{\overset{1}{\wideparen{T+2l_1}}}Z^2(t){\rm d}t\times \int_{\overset{2}{T}}^{\overset{2}{\wideparen{T+2l_2}}}\prod_{r=0}^1Z^2(\vp_1^r(t)){\rm d}t
\end{split}
\ee 
for every fixed $l_1,l_2>0$. 

For simplicity, we use in (\ref{6.3}) an infinite subset of such $l_1^1,l_2^1\in\mbb{R}^+$ that 
\be \label{6.4} 
l_1^1l_2^1=\frac{1}{16\pi}, 
\ee  
and, simultaneously, we make the substitution 
\be \label{6.5} 
T=(x\tau)^{1/\Delta}. 
\ee  
We obtain the next functional. 

\begin{mydef14}
\be \label{6.6} 
\begin{split} 
& \lim_{\tau\to\infty}\frac{1}{\tau}
\left\{
\sum_{\bar{t}_{2\nu+1}\geq (x\tau)^{1/\Delta}}^{\leq (x\tau)^{1/\Delta}+H((x\tau)^{1/\Delta})}Z'(\bar{t}_{2\nu+1})
\right\}\times \\ 
& \left\{
\int_{[(x\tau)^{1/\Delta}]^1}^{[(x\tau)^{1/\Delta}+2l_1^1]^1}Z^2(t){\rm d}t \times 
\int_{[(x\tau)^{1/\Delta}]^2}^{[(x\tau)^{1/\Delta}+2l_2^1]^2}\prod_{r=0}^{1}Z^2(\vp_1^r(t)){\rm d}t
\right\}^{-1}=x
\end{split} 
\ee
for every fixed $x,l_1^1,l_2^1>0$ and admissible $\Delta\in (0,1/6]$.  
\end{mydef14} 

In the special case of Fermat's rationals, see (\ref{4.9}), we have the following. 

\begin{mydef43}
\be \label{6.7} 
\begin{split}
& \lim_{\tau\to\infty}\frac{1}{\tau}
\left\{
\sum_{\bar{t}_{2\nu+1}\geq (\FR\tau)^{1/\Delta}}^{\leq (\FR\tau)^{1/\Delta}+H((\FR\tau)^{1/\Delta})}Z'(\bar{t}_{2\nu+1})
\right\}\times \\ 
& \left\{
\int_{[(\FR\tau)^{1/\Delta}]^1}^{[(\FR\tau)^{1/\Delta}+2l_1^1]^1}Z^2(t){\rm d}t\times 
\int_{[(\FR\tau)^{1/\Delta}]^2}^{[(\FR\tau)^{1/\Delta}+2l_2^1]^2}\prod_{r=0}^{1}Z^2(\vp_1^r(t)){\rm d}t
\right\}^{-1}= \\ 
& \FR
\end{split}
\ee 
for every fixed Fermat's rational, $l_1^1,l_2^1>0$, and admissible $\Delta\in (0,1/6]$. 
\end{mydef43}

Consequently, the next theorem follows from (\ref{6.7}). 

\begin{mydef15}
The $\zeta$-condition 
\be \label{6.8} 
\begin{split}
& \lim_{\tau\to\infty}\frac{1}{\tau}
\left\{
\sum_{\bar{t}_{2\nu+1}\geq (\FR\tau)^{1/\Delta}}^{\leq (\FR\tau)^{1/\Delta}+H((\FR\tau)^{1/\Delta})}Z'(\bar{t}_{2\nu+1})
\right\}\times \\ 
& \left\{
\int_{[(\FR\tau)^{1/\Delta}]^1}^{[(\FR\tau)^{1/\Delta}+2l_1^1]^1}Z^2(t){\rm d}t\times 
\int_{[(\FR\tau)^{1/\Delta}]^2}^{[(\FR\tau)^{1/\Delta}+2l_2^1]^2}\prod_{r=0}^{1}Z^2(\vp_1^r(t)){\rm d}t
\right\}^{-1}\not= \\ 
& 1
\end{split}
\ee 
on the set of all Fermat's rationals and every fixed $l_1^1,l_2^1>0$, and admissible $\Delta\in (0,1/6]$ expresses the next $\zeta$-equivalent of the Fermat-Wiles theorem. 
\end{mydef15} 

\begin{remark}
The results corresponding with the partition $3=1+1+1$ can be obtained as those ones in (\ref{6.1}) -- (\ref{6.8}). 
\end{remark}

I would like to thank Michal Demetrian for his moral support of my study of Jacob's ladders.

\end{document}